\documentclass[10pt,a4paper]{amsart}
\usepackage{amsmath}
\usepackage{amscd}

\newtheorem{thm}{Theorem}[section]
\newtheorem{cor}[thm]{Corollary}
\newtheorem{lem}[thm]{Lemma}
\newtheorem{prop}[thm]{Proposition}

\newcommand{\thmref}[1]{Theorem~\ref{#1}}
\newcommand{\propref}[1]{Proposition~\ref{#1}}

\newsavebox{\SmallMathBox}

\def\hh{{\rm H}}

\def\LIM{{\rm LIM}}

\def\Ci{C^\infty}

\def\dd{\partial}
\def\db{\bar \partial}
\def\Di{D\kern -.65em /}
\def\Dii{D\kern -.45em /}
\def\di{{\dd}\kern -.55em /}
\def\dii{{\dd}\kern -.40em /}

\def\noi{\noindent}

\def\Pf{\mathbb{P}}
\def\wPf{{\widehat \Pf}}
\def\wP{{\widehat P}}
\def\Sf{\mathbb{S}}

\def\wD{{\widehat D}}
\def\bD{{\bf D}}

\def\Df{\mathbb{D}}
\def\Dr{D^{r}}
\def\Dfr{\mathbb{D}^{(r)}}
\def\wDf{{\widehat \Df}}

\def\wDd{{\widehat \Dd}}
\def\Gr{\mbox{{${\mathcal Gr}_{\infty}$}}}
\def\Gra{\mbox{{${\mathcal Gr}_{-1}$}}}

\def\wpGr{\mbox{{${\widehat {\mathcal Gr}}_{-1}$}}}
\def\wsGr{\mbox{{${\widehat {\mathcal Gr}}_{\infty}$}}}

\def\wPdf{\mbox{P($\wDf$)}}

\def\to{\rightarrow}
\def\too{\longrightarrow}

\def\RR{{\bf R}}

\def\ZZ{{\bf Z}}

\def\Cc{{\mathcal C}}
\def\Dd{{\mathcal D}}

\def\Ff{{\mathcal F}}

\def\Hh{{\mathcal H}}

\def\Kk{{\mathcal K}}

\def\Nn{{\mathcal N}}

\def\Rr{{\mathcal R}}
\def\Ss{{\mathcal S}}

\def\Ww{{\mathcal W}}

\def\={\cong}
\def\>{\supset}
\def\<{\subset}
\def\ii{^{-1}}
\def\12{\frac{1}{2}}
\def\2{\Dd}
\def\3{\Nn}
\def\4{\Rr}
\def\6{\cup}
\def\8{\otimes}
\def\0{^{\circ}}

\def\a{\alpha}

\def\th{\theta}

\def\la{\lambda}

\def\Si{\Sigma}

\def\z{\zeta}

\def\Det{\mbox{\rm Det\,}}
\def\DET{\mbox{\rm DET\,}}

\def\dom{\mbox{\rm dom\,}}

\def\End{\mbox{\rm End}}

\def\Hom{\mbox{\rm Hom}}

\def\index{\mbox{\rm index\,}}

\def\Ker{\mbox{\rm Ker}}

\def\ran{\mbox{\rm range\,}}

\def\Si{S\kern -.65em /}

\def\Tr{\mbox{\rm Tr\,}}

\begin{document}

\title[Determinants of EBVPs and the Quillen metric]{Relative Zeta Determinants and the Quillen
Metric}

\author{Simon Scott}

\address{Department of Mathematics, King's College, London WC2R 2LS, U.K.}

\email{sscott@mth.kcl.ac.uk}

\begin{abstract} We compute the relation between the Quillen metric and
and the canonical metric on the determinant line bundle for a
family of elliptic boundary value problems of Dirac-type. To do
this we present a general formula relating the $\z$-determinant
and the canonical determinant for a class of higher-order elliptic
boundary value problems.
\end{abstract}

\maketitle


\section{Introduction}

The determinant line bundle for a family of first-order elliptic
operators over a closed manifold was first studied in remarkable
papers by Quillen \cite{Qu85}, Atiyah and Singer \cite{AtSi84} ,
Bismut \cite{Bi86}, Bismut and Freed \cite{BiFr86}. Using spectral
$\z$-functions it was shown that the regularized geometry of the
determinant line bundle encodes a subtle relation between the
Bismut local family's index theorem and fundamental non-local
spectral invariants. A great deal of subsequent work has been
directed towards a generalization of this theory to manifolds with
boundary. The corresponding object of interest is the determinant
bundle for a family of elliptic boundary value problems (EBVPs).
The choice of boundary conditions for the family introduces a new
degree of complexity into the problem, or, more optimistically, a
new degree of freedom. Indeed, explicit computations of
$\z$-determinants on closed manifolds are generally achieved by an
identification with a spectrally equivalent EBVP. The crucial fact
behind this is a canonical identification of the space of interior
solutions to the Dirac operator with their boundary traces defined
by the {\em Poisson operator}, enabling one to resolve the EBVP by
solving an equivalent problem for pseudodifferential operators on
the space of boundary sections.

In \cite{Sc99} (math/9812124) it was explained how these facts
lead to a regularization procedure for determinants of Dirac type
operators and a corresponding regularized geometry of the
determinant line bundle for families of EBVPs. In this note simple
formulas are given relating the canonical metric and the Quillen
$\z$-function metric. Details of the constructions here along with
the relation between the $\z$-function and canonical curvature
will be presented in \cite{Sc2000}.

Both the canonical regularization and the $\zeta$-function
regularization depend on choices. The former depends on an
admissible choice of a basepoint boundary condition, the
corresponding canonical determinant is essentially the quotient of
the (unregularized) determinant of the EBVP with that defined by
the basepoint condition.  The later depends on a choice of
spectral cut.  Up to the choice of basepoint the construction of
the canonical regularization is purely topological, being achieved
via natural isomorphisms of the determinant line bundle.
Furthermore, it is a completely algebraic operator-theoretic
honest Fredholm determinant. The first fact is encouraging that
this `god-given' regularization may be closely related to the
$\z$-determinant. The second fact is rather less encouraging, the
$\z$-determinant is defined by an analytic regularization and its
meromorphic continuation relies on delicate asymptotic properties.
Both points, however, are correct. The $\z$-determinant really is
a subtle analytic object not accessible via operator determinants,
however the {\em relative} $\z$-determinant with respect to
boundary conditions $P^1,P^2$ is actually not mysterious and
really is given by a Fredholm determinant over the boundary.

To compute the relation between the metrics we have to study the
Dirac Laplacian EBVP, which  {\it a priori} is a quite complicated
operator since one must additionally restrict the range of the
EBVP for the Dirac operator to lie in the domain of the adjoint
EBVP. However, two key ideas reduce this to a computable problem.
First, we use a useful trick which identifies a certain class of
higher-order EBVPs with an equivalent first-order system (see
Grubb \cite{Gr99} for a recent similar method used to compute the
asymptotic expansion of the heat kernel via `doubling-up').
Second, we use an explicit formula for the {\em relative
parametrix} of the Schwartz kernel of an EBVP relative to
different choices of boundary condition which generalizes the
constructions in \cite{ScWo99}. One is therefore reduced to
computing only a {\em relative determinant}, and this
serendipitously depends only on boundary data.

The results in this paper are closely connected with those of
Mueller \cite{Mu98} on relative determinants on non-compact
manifolds. This will be discussed in \cite{Sc2000}.  The general
results presented here are a generalization of the formulas in the
paper \cite{ScWo99} of K.Wojciechowski and the author in which the
special case of first-order self-adjoint EBVPs over an
odd-dimensional manifold was treated via variational methods. Here
we remove those three  conditions:
first-order,self-adjoint,odd-dimensions. I refer to \cite{Sc99}
and \cite{ScWo99} for background material and to \cite{ScTo99} for
a detailed presentation of the 1-dimensional case. I thank Gerd
Grubb for helpful conversations.


\section{The Determinant Line Bundle for Families of Higher-order
EBVPs}

 We consider a
smooth fibration of manifolds $\pi : Z\too B$ with fibre
diffeomorphic to a compact manifold $X$ with boundary $\dd X = Y$,
and endowed possibly with a vertical Hermitian coefficient bundle
$\xi \to Z/B$ with compatible connection, and such that the
tangent bundle along the fibres $T(Z/B)$ is oriented, spin and
endowed with a Riemannian metric $g_{T(Z/B)}$. This data defines a
corresponding family of Dirac operators  $\Df_{0}=\{ D_{b}:b\in B
\}:\Ff^{0}\too\Ff^{1}$. Here $\Ff^{i}$ are the
infinite-dimensional bundles on $B$ with fibre at $b$  the Frechet
space of smooth sections $ C^{\infty}(M_{b},S^{i}_{b})$, where
$S_{b}^{i}$ are the appropriate Clifford bundles. Since the
dimension of $X$ is unrestricted here, we shall not specify
whether $\Df$ is the family of total Dirac operators or, in even
dimensions, a family of chiral Dirac operators acting between the
bundles of positive and negative chirality fields.

The corresponding structures are inherited on the boundary
fibration $\dd \pi : \dd Z\too B$ of closed manifolds with fibre
$Y$.  We assume a collar neighbourhood $V = [0,1]\times \dd Z $ of
$\dd Z$ in which $g_{T(Z/B)}$ has the product form $ du^{2} +
g_{Y/B}$, where $u$ is a normal coordinate to $\dd Z$ and
$g_{Y/B}$  the induced metric on the boundary fibration, and that
the Hermitian metrics, connections on the bundles $S(Z/B), \xi$
split similarly. In $V_{b} = [0,1]\times \dd X_{b} $ the operators
$D_{b}$ have the form $\sigma_{b}(y)(\dd/\dd u +  A_{b})$ where
$\sigma_{b} : S_{Y_{b}} \to S_{Y_{b}}$ is a unitary bundle
isomorphism and $A_{b}$ is an elliptic self-adjoint operator
identified with the Dirac operator over $Y_{b}$. Each fibre of the
bundle $\Hh_{Y}$ is endowed with a $\ZZ_{2}$-grading $\Hh_{Y,b} =
\hh_{b}^{+}\oplus \hh_{b}^{-}$ where $\hh_{b}^{+}$ (resp.
$\hh_{b}^{-}$) is the direct sum of the eigenspaces of $A_{b}$
with non-negative (resp. negative) eigenvalues. Associated to the
gradings we have a (smooth) Grassmann bundle $\Gr (\Hh_{Y})\too B$
with fibre the infinite-dimensional {\it smooth Grassmannian}
$\Gr(\Hh_{Y,b})$ parameterising projections $P\in\End (\hh_{b})$
which differ by a smoothing operator from the projection
$\Pi_{\geq,b}$ with image $\hh_{b}^{+}$, where projection means
self-adjoint indempotent. $\Gr(\Hh_{Y,b})$ is a dense submanifold
of the larger Grassmannian $\Gra(\Hh_{Y,b})$ parameterising
projections such that $P-\Pi_{\geq,b}$ is a pseudodifferential
operator of order -1.  A {\em Grassmann section} for the family
$\Df_{0}$ is defined to be a smooth section $\Pf$ of the fibration
$\Gr(\Hh_{Y})$, and we denote the space of Grassmann sections by
$\Gr^{(0)}(Y/B)$. Such sections always exist.

More generally, for $i=0,\ldots,r-1$, let $\Df_{i}=\{ D_{i,b}:b\in
B \}:\Ff^{i}\to \Ff^{i+1}$ be families of compatible Dirac
operators over $X$, defined as above. The bundle $\Ff^{i}$ has
fibre $\Ff_{b}^{i} = C^{\infty}(M_{b},S^{i}_{b})$, with
$S^{i}_{b}$ the corresponding bundles of Clifford-modules. We then
consider the family of $r^{th}$ order elliptic differential
operators
\begin{equation}\label{e:rth}
  \Dfr = \Df_{r-1}\cdot\Df_{r-2}\cdot\ldots\cdot\Df_{0}
  : \Ff^{0}\too\Ff^{r}.
\end{equation}
Thus $\Dfr$ parameterizes the operators $\bD_{b} = D_{r,b}\cdot
D_{r-1,b}\cdot\ldots\cdot D_{1,b}
  :\Ff_{b}^{0}\too \Ff_{b}^{r} .$
  In
$V_{b} = [0,1]\times \dd X_{b} $ the operators $D_{i,b}$ have the
 product form $\sigma_{b}(y)(\dd/\dd u + A_{i,b})$.
 For odd $r=2k+1$ we assume that $\Ff^{k}=\Ff^{k+1}$.

For each $\Df_{i}$ we have an associated space $\Gr^{(i)}(Z/B)$ of
Grassmann sections, and to an $r$-tuple $\Pf^r =
(\Pf_{0},\ldots,\Pf_{r-1})$ of Grassmann sections $\Pf_{i}\in
\Gr^{(i)}(Z/B)$ we have a family $(\Df^{(r)},\Pf^r)$ of EBVPs
parameterizing the operators
\begin{equation}\label{e:rbvp}
  \bD_{P_{b}^r} = \bD_{b}: \dom(\bD_{P_{b}^{r}})
  \too L^{2}(X_{b};S^{r}_{b}),
\end{equation}
$$\dom(\bD_{P_{b}^r}) = \{s\in H^{r}(X_{b};S_{b}^{0}):
P_{b,r-1}\gamma_{r-1}(D_{b,r-1}\ldots D_{b,1}s) = 0,\ldots ,
P_{b,0}\gamma_0 s = 0\}$$ where $P_b =
(P_{b,0},\ldots,P_{b,r-1})$, and  $\gamma_{i} :
H^{i+1}(X_b;S_{b}^{i})\to H^{i+1/2}(Y_b,S^{i}_{b,Y})$ is the
restriction operator.

To analyze $(\Df^{(r)},\Pf^r)$ we consider the `equivalent' family
of first-order Dirac type operators
\begin{equation}\label{e:Dequiv}
\wDf = \left(\begin{array}{ccccc} 0   &  0 & \ldots  & 0 &
\Df_{r-1}\\ 0 & 0 & \ldots & \Df_{r-2} & -I\\ \vdots  & \vdots  &
\vdots & \vdots & \vdots\\ \Df_0 & -I   &  0 & \ldots  &
0\end{array}\right) : \Ff_0\oplus\Ff_1\oplus\ldots\Ff_{r-1} \too
\Ff_{r}\oplus\Ff_{r-1}\oplus\ldots\Ff_1
\end{equation}

In $V_{b} = [0,1]\times \dd X_{b} $ the operators $\wD_{b}$
parameterized by $\wDf$ are of Dirac-type with the product form
$\widehat{\sigma}_{b}(y)(\dd/\dd u + \widehat{A}_{b})$, where
$\widehat{\sigma}_{b}:\bigoplus_{i=1}^{r}S^{i}_{b}\to
\bigoplus_{j=0}^{r-1}S^{r-j+1}_{b}$ is a unitary bundle
isomorphism, and $\widehat{A}_{b} = (\bigoplus_{i=0}^{r-1}A_{i,b})
+ R$, where $R$ is an operator of order $0$ on the boundary fields
built from the $\sigma_{i,b}(y)$. We therefore have for each
operator $\wD_{b}$ a Grassmannian $\wpGr(\Hh_{Y,b})$ of
projections associated to the spectral grading $\Hh_{Y,b} =
\hh_{b}^{+}\oplus \hh_{b}^{-}$ defined by the elliptic
self-adjoint operator $\bigoplus_{i=0}^{r-1}A_{i,b}$, where
$\Hh_{Y,b} = \bigoplus_{i=0}^{r-1}L^2 (Y_{b};S^{i}_{Y,b})$.
Globally, we obtain a Grassmann bundle, a space of Grassmann
sections $\wpGr(Z/B)$, and for each $\wPf\in\wpGr(Z/B)$ a family
of first-order EBVPs $(\wDf,\wPf)$ and a family of $r^{th}$ order
EBVPs $(\Df^{(r)},\wPf)$ with $\dom(\bD_{P_{b}}) = \{s\in
H^{r}(X_{b};S_{b}^{0}): \wP_{b}\gamma_{\Dr}s = 0\}$, where
$\gamma_{\Dr}s = (\gamma_0 (s), \gamma_1 (D_{b,1} s),
\ldots,\gamma_{r-1}(D_{b,r-1}\ldots D_{b,1}s))$.

A Grassmann section $\Pf\in \wpGr(Z/B)$ is equivalent to a smooth
ungraded Frechet subbundle $\Ww\too B$ of $\Hh_{Y}$, and for each
such pair of Grassmann sections $\Pf^{0},\Pf^{1}$ we have the
smooth family of Fredholm operators $$(\Pf^{0},\Pf^{1})\in
C^{\infty}(B;\Hom (\Ww^{0},\Ww^{1})),\;\;\;
(\Pf^{0},\Pf^{1})_{b}\equiv P^{1}_{b}P^{0}_{b}: W_{b,0}\too
W_{b,1}, $$ where $\Ww^{i}$ are the bundles defined by $\Pf^i$,
and hence a (Segal) determinant line bundle
$\DET(\Pf^{0},\Pf^{1})$ with its canonical determinant section
$b\mapsto \det(P^{1}_{b}P^{0}_{b})$. Associated to $\wDf$ we have
a canonical Grassmann section $\wPdf$, equal at $b$ to the
Calderon projection $P(\wDd_{b})$ with range equal to the space of
boundary traces $K_{b} = K(\wDd_{b}) = \gamma\Ker
(\wDd_{b})\subset H_{Y_{b}}$ of solutions to the Dirac operator,
where $\gamma$ is the restriction operator to the boundary of
$X_{b}$. For each $\Pf\in \wpGr(Y/B)$ the pair $(\wDf,\Pf)$
therefore has a canonically associated Fredholm family
\begin{equation}\label{e:S(P)}
\Sf(\Pf) := (\wPdf,\Pf): K(\wDf) \to \Ww,
\end{equation}
where $K(\wDf)$ has fibre $K_{b}$, with determinant line bundle
$\DET(\Sf(\Pf))$ with determinant section  $b\mapsto
\det(\Ss(P_{b}))$, where $\Ss(P_{b}):= P_b P(\wDd_{b})$. On the
other hand, we have the smooth Fredholm families $(\Dfr,\Pf)$ and
$(\wDf,\Pf)$, with associated determinant line bundles
$\DET(\Dfr,\Pf),\DET(\wDf,\Pf)$ over $B$ endowed with their
respective determinant sections $b\mapsto \det(\bD_{P_{b}})$ and
$b\mapsto \det(\wD_{P_{b}})$

\begin{thm}\label{t:detlines}
For $P_b\in \wpGr(\Hh_{Y,b})$, $\bD_{P_{b}}, \wD_{P_{b}}$ are
Fredholm operators with kernel and cokernel consisting of smooth
sections. One has $$\index (\bD_{P_{b}}) = \index (\wD_{P_{b}}) =
\index \Ss(P_{b}).$$ For $\Pf,\Pf_1,\Pf_2\in \wpGr (Z/B)$ there
are canonical isomorphisms of determinant line bundles
\begin{equation}\label{e:det=}
\DET(\Dfr,\Pf) \cong \DET(\wDf,\Pf) \cong \DET(\Sf(\Pf)),
\end{equation}
preserving the determinant sections $\det(\bD_{P_{b}})
\longleftrightarrow \det(\wD_{P_{b}}) \longleftrightarrow
\det(\Ss(P_{b})),$ and
\begin{equation}\label{e:det=2}
\DET(\Dfr,\Pf_1)\cong \DET(\Dfr,\Pf_2) \otimes
\DET(\Pf^{1},\Pf^{2}).
 \end{equation}
\end{thm}
\bigskip

With the identifications of \thmref{t:detlines} at hand, we obtain
a commutative diagram of canonical isomorphisms

 $$\begin{CD} \DET(\Dfr,\Pf^1)
@>{\simeq}>> \DET(\Dfr,\Pf^2)\otimes\DET (\Pf^{1},\Pf^{2})\\
@VV{\simeq}V @VV{\simeq}V \\ \DET(\Sf(\Pf^1))  @>{\simeq}>>
\DET(\Sf(\Pf^2))\otimes\DET (\Pf^{1},\Pf^{2})
\end{CD}$$

\vskip 1mm

\noi in which, by \eqref{e:det=}, the vertical maps take the
determinant sections to each other, while in the bottom map

$$\det(P_{b}^{1}\Ss(P_{b}^{2}))\longleftarrow
\det(\Ss(P_{b}^{2}))\otimes \det (P_{b}^{1},P_{b}^{2}).$$

If we assume the EBVPs $\bD_{P^{i}_{b}}$ are invertible at $b$,
then all the above operators are invertible, and we obtain {\em
two} non-zero canonical elements $\det(\Ss(P^{1}))$ and
$\det(P^{1}\Ss(P^{2}))$ in the complex {\em line}
$\Det(\Ss(P^1))$, where for brevity we have dropped the $b$
subscript. The first element is identified with
$\det(\bD_{P^{1}})\in \Det(\bD_{P_{1}})$ by the isomorphism
\eqref{e:det=}, while the second maps to an element of
$\Det(\bD_{P_{1}})$ we shall denote by $\det(\bD_{(P^{1},P^2)})$.
The two elements therefore differ by the complex number
$\det_{\Cc(P^2)}\bD_{P^1}$ where $$\det(\bD_{P_{1}}) = {\rm
det}_{\Cc(P^2)}(\bD_{P^1}).\det(\bD_{(P^{1},P^2)}).$$
$\det_{\Cc(P^2)}\bD_{P^1}$ is called the {\it canonical
regularization of the determinant of the $r^{th}$ order EBVP
$\bD_{P^1}$ relative to the basepoint $P^{2}$}. Notice that this
is a purely {\em topological} regularization. We compute:

\begin{lem}
\begin{equation}\label{e:cdet}
{\rm det}_{\Cc(P^2)}\bD_{P^1} = {\rm
det}_{Fr,K}\left(\frac{\Ss(P^1)}{P^{1}\Ss(P^{2})}\right).
\end{equation}
\end{lem}
Here $\det_{Fr,K}$ means the Fredholm determinant taken on the
Calderon subspace $K = K(\wD)$ defined by $\wD$, and the operator
quotient means
$\left(P^{1}\Ss(P^{2})\right)^{-1}\Ss(P^1):K(\wD)\to K(\wD)$.


\section{Computation of the relative zeta-function determinant}

In order to compare the canonical regularization \eqref{e:cdet}
with the $\z$-function regularization we study the relative
$\z$-determinant of the $r^{th}$ order EBVPs
$\bD_{P^{1}},\bD_{P^{2}}$ with $P^{1},P^{2}\in
\wsGr(\Hh_{Y_{b}})$. We assume the operators have a common ray of
minimal growth $l_{\theta} = \{arg\la = \th :\la\in\mathbb{C}\}$.
This means that there is an open neighbourhood of $l_{\th}$
disjoint from the spectrum of the operators $\bD_{P^1}$ and
$\bD_{P^2}$. We further assume that there exists $r>0$ such that
for $|\la |>r$ the $L^{2}$ operator norms along this ray
$\|(\bD_{P^i}-\la)\ii\|$ are $O(1/|\la |)$.  For $ Re (s)
> 0 $ we can therefore define
\begin{equation}\label{e:Ds}
  \bD_{P^i}^{-s} = \frac{1}{2\pi i}
  \int_{\Gamma}\la^{-s}(\bD_{P^i}-\la)\ii\,d\la,
\end{equation}
where $\Gamma$ is the contour beginning at $\infty$ traversing
$l_{\th}$ to a small circle around the origin, anti-clockwise
around the circle, then back along the ray to $\infty$. For
$Re(s)> n/r$, $n = \dim X$, the operator $ \bD_{P^i}^{-s}$ has a
continuous kernel and the {\em spectral $\z$-function}
\begin{equation}\label{e:zeta}
\zeta_{\theta}(s,\bD_{P^i}) = \Tr(\bD_{P^i}^{-s})
\end{equation}
is well-defined and holomorphic.

The analytic continuation of the $\z$-function to the whole of
$\mathbb{C}$ depends (as does the construction of the canonical
regularization) on properties, first, of the {\em Poisson
operator} $\Kk_b = :\Ci(Y_b;\widehat{S}_{b|Y})
  \too \Ci(X_b;\widehat{S}_{b})$, which defines
 an isomorphism $\Kk_b :
K(\wD_b)\to \Ker(\wDd_b)$, and thus a canonical identification of
interior solutions with boundary traces. See
\cite{Se69,BoWo93,ScWo99} for details of the construction. Second,
we need the {\em relative parametrix}, given in
\propref{p:parametrix} below, which tells us how the Schwartz
kernel of the inverse operator changes under a change of boundary
condition. For a trace-class operator
$A:\Ff_{r,b}\oplus\Ff_{r-1,b}\oplus\ldots\Ff_{1,b}\too\Ff_{0,b}\oplus\Ff_{1,b}
\oplus\ldots\Ff_{r-1,b}$ written as an $r\times r$ block matrix of
operators relative to the direct sum decomposition, we use the
notation $[A]_{(1,1)}$ for the component in the top-left $(1,1)$
position. We have:

\begin{prop}\label{p:parametrix}
Let $P^1,P^2\in\wsGr(H_{Y})$. If $\bD_{P^1},\bD_{P^2}$ are
invertible, then
\begin{equation}\label{e:parametrix}
\bD_{P^1}\ii = \bD_{P^2}\ii -
\left[\Kk(P^1)P^1\gamma\wD_{P^2}\ii\right]_{(1,1)}
 \end{equation}
 Let $\la$ be a complex number such that $\bD_{P^1}-\la$ and
 $\bD_{P^2}-\la$ are both invertible and assume that $\Ff_{0,b}=\Ff_{r,b}$.
 Then
 \begin{equation}\label{e:laparametrix}
(\bD_{P^1}-\la)\ii = (\bD_{P^2}-\la)\ii -
\left[\Kk_{\la}(P^1)P^1\gamma\wD_{\la,P^2}\ii\right]_{(1,1)}.
 \end{equation}
 Here $\Kk_{\la}(P) = \Kk_{\la}\Ss_{\la}(P)\ii P$, where $\Kk_{\la}$
  is the Poisson operator of $\wD_{\la}$, $\Ss_{\la}(P) =
  PP(\wD_{\la})$, and
\begin{equation}\label{e:laDequiv}
\wD_{\la} = \left(\begin{array}{ccccc} -\la   &  0 & \ldots  & 0 &
\Df_{r-1}\\ 0 & 0 & \ldots & \Df_{r-2} & -I\\ \vdots  & \vdots  &
\vdots & \vdots & \vdots\\ \Df_0 & -I   &  0 & \ldots  &
0\end{array}\right) :
\Ff_{b,0}\oplus\Ff_{b,1}\oplus\ldots\Ff_{b,r-1} \too
\Ff_{b,0}\oplus\Ff_{b,r-1}\oplus\ldots\Ff_{b,1}.
\end{equation}
The remainder term on the right-side of \eqref{e:parametrix} and
\eqref{e:laparametrix} is a smoothing operator.
\end{prop}
\begin{proof}
We find that $\bD_{P}\ii = \left[\wDd_{P}\ii \right]_{(1,1)}$ and
hence to prove \eqref{e:parametrix} it is enough to prove that
$\wDd_{P^1}\ii = \wDd_{P^2}\ii - \Kk(P^1)P^1\gamma\wD_{P^2}\ii$,
and this follows in a similar way to \cite{ScWo99}. If
$P^1,P^2\in\wsGr(Z/B)$ then $P^{1}(I-P^2)$ is a smoothing
operator, and so the last statement is immediate.
\end{proof}

 In a recent papers of Grubb \cite{Gr99,Gr99b},
generalizing joint work with Seeley \cite{GrSe96}, it was shown
that for $(m+1)r> n = \dim X$ there is as asymptotic expansion as
$\la\to\infty$
\begin{eqnarray}\label{e:asymp}
\Tr (\dd^{m}_{\la}((\bD_{P}-\la)\ii ) & \sim &
\sum_{j=0}^{\infty}(a_j + b_j)(-\la)^{(n-j)/r - m -1} \\ & & +
\sum_{k=0}^{\infty}(c_k\log (-\la) + \tilde{c}_{k})(-\la)^{-k/r -
m -1}, \nonumber
\end{eqnarray}
where the $a_j$ are integrals of densities locally determined by
the symbol of $\bD$ and the $b_j,c_k$ similarly with densities
determined locally by the symbol of $\bD$ and $P$. The
$\tilde{c}_{k}$ are in general globally determined. From this one
obtains the meromorphic continuation of
$\zeta_{\theta}(s,\bD_{P})$ to $ \mathbb{C}$ and the full pole
structure. If the coefficients $c_0, \tilde{c}_{0}$ vanish there
is no pole at $s=0$ and so we can then define the $\z$-
determinant $${\rm det}_{\zeta,\th}(\bD_{P^i}) :=
\exp(-\frac{d}{ds}_{|s=0}\zeta_{\theta}(s,\bD_{P^i})).$$ From
\cite{Gr99, Gr99b} this applies to the case of the Dirac Laplacian
for $P^i\in \wsGr(H_{Y})$ (see Section 4) and also to the
first-order self-adjoint problems considered in \cite{ScWo99}. See
also \cite{Wo99}.

 On the other hand
$\left[\Kk_{\la}(P^1)P^1\gamma\wD_{\la,P^2}\ii\right]_{(1,1)}$ is
a smoothing and hence trace-class operator and its trace it
determined explicitly by only boundary data. Hence the formula
\eqref{e:laparametrix} tells us that the natural object to {\em
compute} is the {\em relative determinant} $${\rm
det}_{\zeta,\th}(\bD_{P_1},\bD_{P_2}) :=
\exp(-\frac{d}{ds}_{|s=0}\zeta_{\theta}(s,\bD_{P_1},\bD_{P_2})),$$
where the {\em relative spectral $\z$-function} $
\zeta_{\theta}(s,\bD_{P^1},\bD_{P^2}) = \Tr(\bD_{P^1}^{-s} -
\bD_{P^2}^{-s})$ is well-defined and holomorphic for $Re (s) > 0$
and has a meromorphic continuation with simple poles at $s= -k/r,
\ k=1,2,\ldots$.

For a function with an  asymptotic expansion
$f(x)_{x\to\infty}\sim \sum_{j=0}^{\infty}a_{jk}x^{\a_j}\log^{k}
x$, where $k=0,1$ and $\infty > \a_0 > \a_1 > \ldots$ and $\a_k\to
-\infty$, the regularized limit $\LIM_{x\to\infty}f(x)$  is
defined to be equal to the constant $a_{00}$ term in the
expansion. We obtain the following fundamental relation between
the $\z$-function regularization and the canonical regularization:

 \vskip 10mm

\begin{thm}\label{t:1} One has:

\begin{equation}\label{e:1}
\frac{{\rm det}_{\zeta,\th}(\bD_{P_1})}{{\rm
det}_{\zeta,\th}(\bD_{P_2})} = {\rm det}_{\Cc(P_2)}(\bD_{P_1})\, .
\,e^{-\LIM_{\la\to\infty}\log{\rm det}_{\Cc(P_2)}(\bD_{\la,P_1})},
\end{equation}
\vskip 6mm
 \noi where ${\rm det}_{\Cc(P_2)}(\bD_{\la,P_1}) =
{\rm det}_{Fr,K}(( P^{1}\Ss_{\la}(P^{2}))^{-1}\Ss_{\la}(P^1) )$.
If $\zeta_{\theta}(0,\bD_{P^1}) = \zeta_{\theta}(0,\bD_{P^2})$,
then $\LIM$ can be replaced by (the unregularized) $\lim$.
\end{thm}
 \vskip 8mm

\begin{proof} We have that ${\rm det}_{\zeta,\th}(\bD_{P_1},\bD_{P_2}) =
\det_{\zeta,\th}(\bD_{P_1})/\det_{\zeta,\th}(\bD_{P_2})$ while a
direct computation yields
\begin{equation}\label{e:relDs}
  \Tr\left\{(\bD_{P^1}-\la)\ii - (\bD_{P^2}-\la)\ii \right\}
  = \frac{\dd}{\dd \la}\log {\rm
det}_{Fr,K}\left(\frac{\Ss_{\la}(P^1)}{P^{1}\Ss_{\la}(P^{2})}\right)
\end{equation}
from which the equality follows.
\end{proof}

\noi Formula \eqref{e:1} may be regarded as a formula for
determinants in the spirit of the Atiyah-Singer index formula for
elliptic operators: the left-side of \eqref{e:1} is the (relative)
analytical-determinant constructed from $\z$-function analysis,
while ${\rm det}_{\Cc(P_2)}(\bD_{P_1})$ is the (relative)
topological-determinant constructed purely from natural
topological properties of the determinant line bundle.


\section{Applications to the geometry of $\DET(\Df,\Pf)$}

Let $(\Df,\Pf)$ be a family of {\em first-order} EBVPs with
$\Pf\in\Gr(Z/B)$. There is a metric and connection on the
determinant line bundle $\DET(\Df,\Pf)$ associated to the
canonical regularization. Let $\Delta_{P}$ be the Dirac Laplacian
$\Delta = \Dd^{*}\Dd$ with domain $$\dom(\Delta_{P}) = \{s\in
H^{2}(X,S)\to L^{2}(X,S): P\gamma s = 0,\;\;P^{*}\gamma \Dd s =
0\},$$ where $P^{*} = \sigma(I-P)\sigma^{*}$ is the adjoint
boundary condition for $\Dd^{*}$ (not to be confused with the
adjoint of the (self-adjoint) operator $P$) and $S$ the Clifford
bundle of (chiral) spinors. Over the open $U$ subset of $B$ where
the EBVPs are invertible the canonical metric is defined by
\begin{equation}\label{e:cmetric}
  \|{\rm det}(\Dd_{P})\|^{2}_{\Cc} = {\rm det}_{\Cc}(\Delta_{P})
  := {\rm det}_{Fr,K}(\Ss(P)^{*}\Ss(P)),
\end{equation}
where $\Ss(P) := P P(\Dd):K(\Dd)\to \ran(P)$ and $P(\Dd)$ is the
Calderon projection of the (first-order) operator $\Dd$. On the
other hand, $\DET(\Df,\Pf)$ has a Quillen metric, defined by
$\z$-function regularization \cite{Qu85,BiFr86} over $U$ by
\begin{equation}\label{e:zmetric}
  \|{\rm det}(\Dd_{P})\|^{2}_{\z} = {\rm det}_{\z}(\Delta_{P})
  :=  e^{-\zeta_{\Delta_{P}}^{'}(0)}
\end{equation}
where $\zeta_{\Delta_{P}}(s) = \Tr (\Delta_{P}^{-s})$  is defined
around $0$ by analytic continuation.  For the global construction
of these objects see \cite{Qu85,BiFr86,Sc99}. From \thmref{t:1} we
compute:
\bigskip
\begin{thm}\label{t:4}
Over $U\subset B$
\begin{equation}\label{e:metrics}
\frac{\|{\rm det}(\Dd_{P_1})\|_{\z}}{\|{\rm
det}(\Dd_{P_2})\|_{\z}} = \frac{\|{\rm
det}(\Dd_{P_1})\|_{\Cc}}{\|{\rm det}(\Dd_{P_2})\|_{\Cc}} \,.
\end{equation}
\end{thm}
\bigskip
\vskip 5mm

A holomorphic line bundle $L$ endowed with a metric $g$ and local
holomorphic section $s$, has a canonical compatible connection
with curvature 2-form $\RR^{g} = \db\dd\log\|s\|_{g}^{2}$. If we
assume that $\Hh_Y$ is a trivial bundle and that $\Df$ depends
holomorphicaly on $b$, then for two choices of constant Grassmann
sections $\Pf^1 := P^1, \Pf^2 := P^2$ the determinant line bundles
$\DET(\Df,P^i)$ the corresponding $\z$ and $\Cc$ metrics have
respective canonical curvature 2-forms
$\RR_{\z}^{i},\RR_{\Cc}^{i}\in\Omega^{2}(B)$ and we have:

\begin{cor}
With the above assumptions,
\begin{equation}\label{e:curvatures}
\RR^{1}_{\z} - \RR^{2}_{\z} = \RR^{1}_{\Cc} - \RR^{2}_{\Cc}.
\end{equation}
\end{cor}
\noi This would apply, for example, to a family of $\db$-operators
coupled to a Hermitian vector bundle over a Riemann surface with
boundary. The extension of \eqref{e:curvatures} to more general
families of EBVPs will be presented in \cite{Sc2000}.

Finally, we remark the $\z$-determinant formulae recently
published in \cite{LeTo98} for a $r^{th}$ order EBVP in dimension
one follow from \thmref{t:1}, while for self-adjoint first-order
EBVPs over odd-dimensional manifolds (here the condition
$\zeta_{\theta}(0,D_{P^1}) = \zeta_{\theta}(0,D_{P^2})$ applies)
\thmref{t:1} reproves the result of \cite{ScWo99}.


\begin{thebibliography}{99}

\bibitem{AtSi84} Atiyah, M.F., Singer, I.M.: 1984 'Dirac operators coupled to vector
potentials', {\it Proc. Natl. Acad. Sci. {\bf 81}}, 2597--2600.

\bibitem{Bi86} Bismut, J.M.: 1986,
'The Atiyah-Singer index theorem for families of Dirac operators:
Two heat equation proofs', {\it Invent. Math. {\bf 83}}, 91--151.

\bibitem{BiFr86} Bismut, J.M., Freed,D.: 1986 'The analysis of elliptic families:(I)
Metrics and connections on determinant bundles', {\it Commun.
Math. Phys. {\bf 106}}, 159--176.

\bibitem{BoWo93} Boo\ss--Bavnbek, B., and Wojciechowski, K.P.:
1993, {\it Elliptic Boundary Problems for Dirac Operators},
Birkh\"auser, Boston.

\bibitem{Gr99} Grubb, G.: 1999, `Trace expansions for
pseudodifferential boundary problems for Dirac-type operators and
more general systems', {\it Ark. Mat. {\bf 37}}, 45--86.

\bibitem{Gr99b} Grubb, G.: 1999, `Poles of zeta and eta functions
for perturbations of the Atiyah-Patodi-Singer problem', preprint.

\bibitem{GrSe96} Grubb, G., and Seeley. R.: 1996,
`Zeta and eta functions for Atiyah-Patodi-Singer operators', {\it
Invent. Math. {\bf 121}}, 481--529.

\bibitem{LeTo98} Lesch, M., and Tolksdorf. J.: 1998,
`On the determinant of one dimensional elliptic boundary value
problems', {\it Commun. Math. Phys. {\bf 193}}, 643--660.

\bibitem{Mu98} M$\ddot {\rm u}$ller, W.: 1998, ` Relative zeta functions,
relative determinants and scattering theory',{\it Comm. Math.
Phys. {\bf 192}}, 309--347.

\bibitem{Qu85} Quillen, D.G.: 1985, `Determinants of Cauchy-Riemann operators over a
Riemann surface', {\it Funk. Anal. i ego Prilozhenya {\bf 19}},
37--41.

\bibitem{Sc95} Scott, S.G.: 1995, `Determinants of Dirac boundary value
problems over odd--dimensional manifolds', {\it Comm. Math. Phys.
{\bf 173}}, 43--76.

\bibitem{Sc99} Scott, S.G.: 1999, `Splitting the curvature of the determinant line bundle',
{\it Proc. Am. Math. Soc.}, to appear; math/9812124.

\bibitem{Sc2000} Scott, S.G.: In preparation.

\bibitem{ScTo99} Scott, S.G., and Torres, F.: 1999, `Geometry of the determinant
line bundle in dimension one', preprint.

\bibitem{ScWo99} Scott, S.G., and Wojciechowski, K.P.:
[1] 1999, `The $\zeta$--Determinant and Quillen's determinant on
the Grassmannian of elliptic self-adjoint boundary conditions',
{\it C.R. Acad. Sci. Paris, {\bf t.328, Serie I}}, 139-144: [2]
1999, `The $\zeta$--Determinant and Quillen's determinant for a
Dirac operator on a manifold with boundary', {\it GAFA}, to
appear.

\bibitem{Se69} Seeley, R.T.: 1969, `Topics in pseudo-differential operators',
 in {\it Pseudo-Differential Operators (C.I.M.E.,Stresa, 1968) (Nirenburg, L., ed.)
 },43--76, Edizioni Cremonese, Rome.

\bibitem{Wo99} Wojciechowski, K.P.: 1999, `The $\z$-determinant and the additivity
of the $\eta$-invariant on the smooth, self-adjoint
Grassmannian',{\it Comm. Math. Phys. {\bf 169}}, 423--444.

\end{thebibliography}
\end{document}